\documentstyle[11pt,epsf]{article}
\title{On the spectrum of a partial theta function}
\author{Vladimir Petrov Kostov\\ 
Universit\'e de Nice, 
Laboratoire de Math\'ematiques, Parc Valrose,\\ 06108 Nice Cedex 2, France,  
e-mail: kostov@math.unice.fr} 
\date{}
\newtheorem{tm}{Theorem}[section]

\newtheorem{rem}[tm]{Remark}

\newtheorem{lm}[tm]{Lemma}

\newtheorem{nota}[tm]{Notation}
\begin{document} 
\maketitle 
\begin{abstract}
The bivariate series $\theta (q,x):=\sum _{j=0}^{\infty}q^{j(j+1)/2}x^j$ 
defines  
a {\em partial theta function}. For fixed $q$, $\theta (q,.)$ is an 
entire function. We show that for $|q|\leq 0.108$ the function 
$\theta (q,.)$ has no multiple zeros.\\ 


\end{abstract}

\section{Introduction}

Consider the bivariate series $\theta (q,x):=\sum _{j=0}^{\infty}q^{j(j+1)/2}x^j$ 
for $(q,x)\in {\bf C}^2$, $|q|<1$. We consider $x$ as a variable and $q$ 
as a parameter. For each $q$ fixed the series defines an entire function 
called a {\em partial theta function}. This terminology stems from the fact 
that $\theta (q^2,x/q)=\sum _{j=0}^{\infty}q^{j^2}x^j$ while the series 
$\sum _{j=-\infty}^{\infty}q^{j^2}x^j$ defines the Jacobi theta function; in the 
series for $\theta$ only a partial summation (i.e. excluding negative indices) 
is performed. 

The function $\theta$ has been applied in several domains: 
in asymptotic analysis (see \cite{BeKi}), in the theory 
of (mock) modular forms (see \cite{BrFoRh}), in Ramanujan type $q$-series 
(see \cite{Wa}), in statistical physics 
and combinatorics (see \cite{So}), in questions concerning hyperbolic 
polynomials (i.e. real polynomials with all roots real, see \cite{KaLoVi}, 
\cite{KoSh} and \cite{Ko2}). For more information about $\theta$,  
see also~\cite{AnBe}.
 
For real $q$ and $x$ there are countably many values 
$0.3092493386\ldots =\tilde{q}_1<\tilde{q}_2<\cdots <1$ of $q$ such that 
$\theta (q, \cdot )$ has a unique multiple zero of multiplicity $2$. 
Moreover, $\lim _{j\rightarrow \infty }\tilde{q}_j=1$, 
see \cite{KoSh} and \cite{Ko2}. 
These values of $q$ are said 
to belong to the {\em spectrum} of $\theta$. The double zero is the rightmost 
of the real zeros of $\theta$.
For $q\in (\tilde{q}_j,\tilde{q}_{j+1}]$ the 
function $\theta (q, \cdot )$ has exactly $j$ complex conjugate pairs of zeros 
(counted with multiplicity). (It is not clear yet whether 
they are all simple or not.) It is proved in \cite{Ko3} that 
$\tilde{q}_j=1-(\pi /2j)+o(1/j)$ and that the double real zeros of 
$\theta (\tilde{q}_j,.)$ tend to $-e^{\pi}=-23.1407\ldots$ 
(this number appears in a different context, but always in relationship 
with $\theta$, in a theorem 
announced in~\cite{Ka}). 
   
In the present paper we show that when $q$ and $x$ are complex,  
there exists a neighbourhood of $0$ free of spectral points. 
More precisely, we 
prove the following theorem:

\begin{tm}\label{tm1}
There is no spectral value of $q$ for $q\in {\bf C}$, $|q|\leq 0.108$. 
\end{tm}

Our interest in a neighbourhood of $0$ is explained by the fact that 
$\theta (0, \cdot )\equiv 1$ while for any nonzero $q$ on the open unit disk, 
$\theta (q, \cdot )$ 
is a nontrivial entire function. The number $0.108$ certainly does not give 
the best possible estimate for the neighbourhood free from spectral values. 
On the other hand, it cannot be replaced by a number 
greater than or equal 
to $0.3092493386\ldots$ because the latter 
belongs to the spectrum of $\theta$.   

{\bf Acknowledgement.} It seems that A. Sokal was the first to 
ask the question about the existence of a neighbourhood of $0$ 
free of spectral values. On his homepage (see  
http://www.maths. qmul.ac.uk/~pjc/csgnotes/sokal/) one can find conjectures 
about the partial theta function. His claim that the author's result 
should hold for $0.2078750206$ in the place of $0.108$ is supported by 
a sketch of proof (a proof of this is reported to have been given also by 
Jens Forsg{\aa}rd, a student of B.Z. Shapiro).   
The author is grateful to B.Z. Shapiro for the formulation of the problem 
and comments 
on this text, and also to the anonymous referee for his useful remarks. 

\section{Proof of Theorem~\protect\ref{tm1}}

We are looking for a function $\theta$ representable 
in the form of an infinite product 

\begin{equation}\label{product}
\theta =\sum _{s=0}^{\infty}q^{s(s+1)/2}x^s=\prod _{j=1}^{\infty}(1+x/\xi _j)~,
\end{equation} 
where $\{ -\xi _j\}$ is the set of zeros of $\theta (q, \cdot )$. 
(For $q\in (0,0.3092493386\ldots )$ such a presentation exists because 
the function $\theta (q,.)$ is an entire 
function of order $0$ of the Laguerre-P\'olya class $\cal{LP-I}$, see 
\cite{KaLoVi}.) 
For $|q|$ small enough we look for zeros of the form 
$-\xi _j=-1/q^j\Delta _j$ (i.e. $1/\xi _j=q^j\Delta _j$), 
where $\Delta _j$ are complex numbers close to $1$. In other words, the 
zeros of $\theta$ are close to the terms of a geometric progression. When 
$|q|$ is small enough and all $\Delta _j$ are uniformly close to $1$, then 
all zeros are distinct and the corresponding values of $q$ 
are not from the spectrum of $\theta$. Thus we deduce  
Theorem~\ref{tm1} from the following theorem whose proof follows:

\begin{tm}
For $|q|\leq 0.108$ the function $\theta$ is representable in the form 
(\ref{product}) with $1/\xi _j=q^j\Delta _j$, 
where $\Delta _j\in [0.2118,~1.7882]$ 
for $j=1,2,\ldots$.
\end{tm}

The fact that if $q$ and $\Delta _j$ satisfy the conditions of the theorem,  
then all zeros of $\theta$ are distinct, is proved at the end of 
Section~\protect\ref{convergence}.

\section{Formal solution}

We show first that every quantity $\Delta _j$ can be represented as a formal 
power series in $q$. To this end we observe that 
expanding the infinite product in  
(\ref{product}) as a power series in $x$ gives 

$$\begin{array}{ccl}\sum _{s=0}^{\infty}q^{s(s+1)/2}x^s&=&
\sum _{s=0}^{\infty}e_s(1/\xi _1,1/\xi _2,\ldots )x^s\\ \\ 
&=&\sum _{s=0}^{\infty}e_s(q\Delta _1,q^2\Delta _2,\ldots )x^s~,
\end{array}$$
where $e_s$ is the $s$th elementary symmetric function. 
Hence 

\begin{equation}\label{equ1}
e_s(\Delta _1,q\Delta _2,\ldots )=q^{s(s-1)/2}~~,~~s=1,2,\ldots ~.
\end{equation} 
We are going to show that 
the infinite system of these equations can be solved for 
$\Delta _1$, $\Delta _2$, $\ldots$. 
Every equation (\ref{equ1}) is a formal power series (FPS) in the 
infinitely many variables $q$, $\Delta _1$, $\Delta _2$, $\ldots$. 
Set $\Delta ^j:=(\Delta _j,\Delta _{j+1},\ldots )$, $\Delta :=\Delta ^1$. 
We denote by the letter $T$ (indexed or not) an FPS 
in the indicated variables  
$q$ and $\Delta ^j$.  
After division by $q$ the $s=1$ case of equations (\ref{equ1}) reads

\begin{equation}\label{equ2}
1=\Delta _1+q\Delta _2+q^2\Delta _3+\cdots =\Delta _1+qT_{1,2}(q,\Delta ^2)
\end{equation}
The double index of $T_{1,2}$ means that this is an FPS connected with 
the $s=1$ case of equations (\ref{equ1}) 
and in which the 
first of the variables $\Delta _j$ that contributes is $\Delta _2$. 
Similarly, after division by 
$q^{s(s+1)/2}$ equation (\ref{equ1}) for general $s$ becomes

\begin{equation}\label{equ3}
1=\Delta _1\cdots \Delta _s+qT_{s,1}(q,\Delta )
\end{equation}
In what follows we refer to equation (\ref{equ3}) also as to equation 
$(F_s)$. We denote by $(E_{1,2})$ equation (\ref{equ2}) (written in the form 
$\Delta _1=1-qT_{1,2}(q,\Delta ^2)$) and by $(E_{2,2})$ the equation 
obtained from equation $(F_2)$ after substituting in it 
$1-qT_{1,2}(q,\Delta ^2)$
for $\Delta _1$. Hence equation 
$(E_{2,2})$ reads 

$$(E_{2,2})~:~1=(1-qT_{1,2}(q,\Delta ^2))\Delta _2+
qT_{2,1}(q,1-qT_{1,2}(q,\Delta ^2),\Delta ^2)~.$$
Its right-hand side is of the form $\Delta _2+qT^*(q,\Delta ^2)$. 
Hence this can be solved for $\Delta _2$. Indeed, one can  
apply the implicit function theorem here (at $q=0$)  
and obtain an equation of the form

$$(E_{2,3})~:~\Delta _2=1-qT_{2,3}(q,\Delta ^3)~.$$
If we solve $(E_{2,3})$ for $\Delta _2$ and substitute this in $(E_{1,2})$, 
we find the equation
$$(E_{1,3})~:~\Delta _1:=1-qT_{1,3}(q,\Delta ^3)~.$$
In what follows we denote by $(E_{s,r})$ an equation of the form 
$\Delta _s=1-qT_{s,r}(q,\Delta ^r)$. 

The remainder of the proof proceeds by induction. 
Suppose that equations 
$(E_{j,i})$ are constructed for $j=1,\ldots ,s$, $i=j,j+1,\ldots ,s+1$. 
(For $s=2$ we already constructed equations $(E_{1,1})$, 
$(E_{1,2})$, $(E_{1,3})$, $(E_{2,2})$, $(E_{2,3})$.) 

Consider equation $(F_{s+1})$. 
We solve the system of equations $(E_{j,s+1})$ for the variables $\Delta _j$, 
$j=1,\ldots ,s$, and substitute 
this in $(F_{s+1})$. This yields an equation of the form 

$$(E_{s+1,s+1})~:~1=\tilde{\Delta}_1\cdots \tilde{\Delta}_s\Delta _{s+1}+
qT_{s+1,1}(q,\tilde{\Delta}_1,\ldots ,\tilde{\Delta}_s,\Delta ^{s+1})~,$$
where $\tilde{\Delta}_i=1-qT_{i,s+1}(q,\Delta ^{s+1})$. One can express 
$\Delta _{s+1}$ from equation $(E_{s+1,s+1})$ 
(the implicit function theorem is applicable at $q=0$) which gives the equation 
$(E_{s+1,s+2})$. Express then $\Delta _{s+1}$ from it (i.e. 
set $\Delta _{s+1}=1-qT_{s+1,s+2}(q,\Delta ^{s+2})$) and substitute 
$1-qT_{s+1,s+2}(q,\Delta ^{s+2})$ for 
$\Delta _{s+1}$ in equations $(E_{j,s+1})$, $j=1,\ldots ,s$. This gives the 
equations $(E_{j,s+2})$, $j=1,\ldots ,s$.

Applying this above procedure infinitely many times we obtain the equations 
$(E_{s,\infty })$ which express the quantities $\Delta _s$ as FPS 
in $q$ of the form $\Delta _s=1+O(q)$. These FPS stabilize because at every 
substitution of $\Delta _s$ by $1-qT_{s,r}(q,\Delta ^r)$ 
the power of $q$ increases.

\begin{rem}
{\rm It is easy to deduce from the above reasoning 
that all coefficients of $\Delta _s$ 
(when considered as power series in $q$) are integer. 
We list below the first 10 coefficients of 
$\Delta _1$, $\ldots$, $\Delta _5$:}

$$\begin{array}{rrrrrrrrrr}
1&-1&-1&-1&-2&-4&-10&-25&-66&-178\\ 1&0&0&1&3&9&24&66&180&498\\ 
1&0&0&0&0&0&-1&-3&-9&-22\\ 1&0&0&0&0&0&0&0&0&0\\ 
1&0&0&0&0&0&0&0&0&0\end{array}$$ 
{\rm  
It would be interesting to (dis)prove that 
$\Delta _s=1+(-1)^sq^{\kappa _s}\Phi _s$, where $\Phi _s$ is an FPS with 
positive coefficients (for $s=1$ this is proved in \cite{So}), 
and the natural numbers $\kappa _s$ 
form an increasing sequence. It is clear that $\kappa _1=1$, $\kappa _2=3$ and 
$\kappa _3=6$. It would be interesting to explicit $\kappa _s$. 
A combinatorial interpretation of the coefficients of $\Delta _1$ is given 
in~\cite{Pr}.}
\end{rem}

\section{Proof of the convergence\protect\label{convergence}}

\begin{nota}
{\rm We denote by $U$ the infinite column vector whose entries equal $1$ 
(i.e. $U=(1,~1,~\ldots ~)^T$) and similarly we set 
$V:=(1,~\Delta _1,~\Delta _1\Delta _2,
~\Delta _1\Delta _2\Delta _3,~\ldots ~)^T$. 
Denote by $\sigma _s$ the right-hand side of equation (\ref{equ1}) divided 
by $q^{s(s+1)/2}$ (hence $\sigma _s=\Delta _1\Delta _2\ldots \Delta _s+
q\Delta _1\Delta _2\ldots \Delta _{s-1}\Delta _{s+1}+O(q^2)$) 
and by $L_s$ the infinite 
square matrix with $1$ on the diagonal, with 
$q^{s-1}\Delta _s$, $q^{s-2}\Delta _s$, $\ldots$, $q\Delta _s$, $0$, $0$, $\ldots$ 
on the first subdiagional and with zeros elsewhere. That is, 
$L_1=I=$diag$(1,1,\ldots )$,}

$$L_2=\left( \begin{array}{cccc}1&0&0&\cdots\\q\Delta _2&1&0&\cdots\\
0&0&1&\cdots\\\vdots&\vdots&\vdots&\ddots\end{array}\right) ~~,~~
L_3=\left( \begin{array}{ccccc}1&0&0&0&\cdots\\q^2\Delta _3&1&0&0&\cdots\\
0&q\Delta _3&1&0&\cdots\\0&0&0&1&\cdots\\
\vdots&\vdots&\vdots&\vdots&\ddots\end{array}\right) ~~{\rm etc.~~Hence}$$
\end{nota}
$$\begin{array}{rcl}L_2V&=&(1,~\Delta _1+q\Delta _2,~\Delta _1\Delta _2,~ 
\Delta _1\Delta _2\Delta _3,
~\Delta _1\Delta _2\Delta _3\Delta _4,~\ldots ~)^T~,\\ \\  
L_3L_2V&=&(1,~\Delta _1+q\Delta _2+q^2\Delta _3,~\Delta _1\Delta _2+
q\Delta _1\Delta _3+q^2\Delta _2\Delta _3,~\Delta _1\Delta _2\Delta _3,~ 
\Delta _1\Delta _2\Delta _3\Delta _4,~\ldots ~)^T\end{array}$$
and so on. It is easy to see that 

$$\cdots L_4L_3L_2V=(1,~\sigma _1,~\sigma _2,~\sigma _3,~\ldots ~)^T~.$$
Indeed, if $\tilde{\sigma}_j^k$ denotes the $j$th elementary symmetric 
polynomial of the quantities $\Delta _1$, $q\Delta _2$, $\ldots$, 
$q^{k-1}\Delta _k$, then for $j\leq k$ 
$$\tilde{\sigma}_j^k=\tilde{\sigma}_j^{k-1}+
q^{k-1}\Delta _k\tilde{\sigma}_{j-1}^{k-1}~~~{\rm and}~~~
\tilde{\sigma}_j^k/q^{j(j-1)/2}=\tilde{\sigma}_j^{k-1}/q^{j(j-1)/2}+
q^{k-j}\Delta _k(\tilde{\sigma}_{j-1}^{k-1}/q^{(j-1)(j-2)/2})~.$$
Thus the $(j+1)$st component of the vector $L_{k}\cdots L_2V$ equals 
$\tilde{\sigma}_j^k/q^{j(j-1)/2}$. This is a polynomial in the variables 
$q$, $\Delta _1$, $\ldots$, $\Delta _k$. As $k\rightarrow \infty$, 
it stabilizes 
as a formal power series in the infinitely many variables 
$q$, $\Delta$ and tends to 
$\sigma _j$. (Stabilization is due to the increasing powers of $q$.)  
Hence the system of equations $(F_s)$, $s=0,1,\ldots$ (we set $(F_0):1=1$) 
reads  

\begin{equation}\label{VLU}
\cdots L_4L_3L_2V=U~~,~~{\rm i.e.}~~V=L_2^{-1}L_3^{-1}L_4^{-1}\cdots U~.
\end{equation}
We represent the matrix $L_s$ in the form $L_s=I+N_s$ (where $(N_s)^s=0$). 
In particular, $N_1=0$, 
$N_2=\left( \begin{array}{cccc}0&0&0&\cdots\\q\Delta _2&0&0&\cdots\\
0&0&0&\cdots\\\vdots&\vdots&\vdots&\ddots\end{array}\right)$ etc. 
Hence $L_s^{-1}=I+\sum _{j=1}^{s-1}(-N_s)^j$. 
The following lemma is proved in the next section:

\begin{lm}\label{lemma1}
The entry $(L_s^{-1})_{\mu ,\nu}$ of the matrix $(L_s)^{-1}$
equals

$$\left\{ \begin{array}{cl}
(-1)^{\mu -\nu}\Delta _s^{\mu -\nu}q^{(\mu -\nu )(s-\mu +1)+
(\mu -\nu )(\mu -\nu -1)/2}& 
{\rm for~~}\nu \leq \mu \leq s,\\ 
0&{\rm otherwise.}
\end{array}\right.$$
\end{lm}

We are going now to justify the convergence of the formal series in $q$ 
expressing the quantities $\Delta _j$. Denote by $\beta \in (0,0.7882]$ 
a number such that $|\Delta _j|\in [1-\beta ,1+\beta ]$, $j=1, 2, \ldots$. 
Set $u:=1+\beta$. 
Assume that $|q|\leq a$, $a\in (0,0.108]$. Under these assumptions we give an 
estimation of the moduli of the entries of the matrix 
$L:=L_2^{-1}L_3^{-1}L_4^{-1}\cdots$ that are below the main diagonal. 
The above lemma 
implies our next lemma:

\begin{lm}\label{lemma2}
For $\nu <\mu \leq s$ one has 
$|L_{\mu ,\nu}|\leq u^{\mu -\nu}a^{(\mu -\nu )(s-\mu +1)+
(\mu -\nu )(\mu -\nu -1)/2}$. 
\end{lm}

To obtain a majoration for the entries of $L$ in its $s$th row one can 
\vspace{1mm}

1) ignore the presence of the factors $L_2^{-1}$, $\ldots$, $L_{s-1}^{-1}$ 
(because their $s$th rows have just $1$ in position $s$ and zeros elsewhere) 
and  

2) ignore the rows of the matrices $L_{s+1}^{-1}$, $L_{s+2}^{-1}$, $\ldots$ below 
the $s$th one. 
\vspace{1mm}
Therefore in what follows, instead of $L$ we consider  
the $s\times s$-matrix 
$K:=\tilde{L}_s\tilde{L}_{s+1}\tilde{L}_{s+2}\cdots$, where $\tilde{L}_j$
is the left upper $s\times s$-minor of the matrix $L_j^{-1}$ ($j\geq s$). 
 
Denote by $M$ an $s\times s$-matrix having on its first subdiagonal the entries 
$a^{s-1}u$, $a^{s-2}u$, $\ldots$, $au$ and zeros elsewhere. 
It is clear that $M^s=0$ and that the nonzero entries of $M$ 
are majorations of the moduli 
of the respective entries of the left upper $s\times s$-minor of the 
matrix $N_s$. Hence the moduli of the 
entries of the matrix 
$\tilde{L}_s$ are majorized by the entries of the matrix 
$I+\sum _{j=1}^{s-1}M^j$. In the same way the moduli of the 
entries of the matrix 
$\tilde{L}_k$, $k\geq s$, are majorized by the entries of the matrix 
$I+\sum _{j=1}^{s-1}(a^{k-s}M)^j$. Therefore the moduli of the entries of the 
$s$th row of $K$ (hence of $L$ as well) 
are majorized by the entries of the $s$th row of the product 

$$\Pi (a, M):=\prod _{k=s}^{\infty}\left( I+\sum _{j=1}^{s-1}(a^{k-s}M)^j\right)$$
which (taking into account that $M^s=0$) we represent in the form 
$I+b_1M+b_2M^2+\cdots +b_{s-1}M^{s-1}$. 
 
\begin{lm}\label{lemma3} 
One has $b_j\leq 1/(1-a)(1-a^2)\cdots (1-a^j)$, $j=1,\ldots ,s-1$.
\end{lm}

The lemma is proved in the next section. 
The Lemmas~\ref{lemma1}--\ref{lemma3} 
imply the inequality (where $1\leq \nu \leq s-1$)  

$$|L_{s,\nu}|\leq 
\frac{u^{s-\nu}a^{(s-\nu )+(s-\nu )(s-\nu -1)/2}}{(1-a)(1-a^2)\cdots (1-a^{s-\nu})}=
\frac{u^{s-\nu}a^{(s-\nu )(s-\nu +1)/2}}{(1-a)(1-a^2)\cdots (1-a^{s-\nu})}~.$$
Hence the equation (\ref{VLU}) for $V$ implies that 

$$1-\sum _{\nu =1}^{s-1}|L_{s,\nu}|\leq 
|\Delta _1\cdots \Delta _s|\leq 1+\sum _{\nu =1}^{s-1}|L_{s,\nu}|~.$$
The following two inequalities (resulting from the conditions $a\in (0,0.108]$ 
and $\beta \in (0,0.7882]$) will be used in our estimates:

\begin{equation}\label{twoconditions}
0<a<1/3~~~~{\rm and}~~~~0<au<1~.
\end{equation}
Hence for $\nu =s-1$ (resp. $\nu =s-2$) one has 
$|L_{s,s-1}|\leq ua/(1-a)$ 
(resp. 
$|L_{s,s-2}|=u^2a^2(a/(1-a)(1-a^2))
\leq u^2a^2((1/3)/(2/3)(8/9))<u^2a^2$).  
For $\nu \leq s-3$ it is true that 

$$\frac{a^{(s-\nu )(s-\nu -1)/2}}{(1-a)(1-a^2)\cdots (1-a^{s-\nu})}=
a^{(s-\nu )(s-\nu -3)/2}\left( \frac{a}{1-a}\right) \left( 
\frac{a}{1-a^2}\right) \cdots \left( \frac{a}{1-a^{s-\nu}}\right) <1$$
because $a\in (0,1)$ (hence $a^{(s-\nu )(s-\nu -3)/2}\in (0,1)$) 
and all fractions belong to $(0,1)$ 
due to $a\in (0,1/3)$.  
Hence for $\nu \leq s-2$ the inequalities 
$$|L_{s,\nu}|\leq (ua)^{s-\nu}~~{\rm and}~~ 
\sum _{\nu =1}^{s-2}|L_{s,\nu}|\leq (ua)^2+(ua)^3+\cdots =(ua)^2/(1-ua)$$ 
hold true. Thus 
$$1-ua/(1-a)-(ua)^2/(1-ua)\leq 
|\Delta _1\cdots \Delta _s|\leq 1+ua/(1-a)+(ua)^2/(1-ua)~.$$ 
Observe that the right-hand and left-hand sides do not depend on $s$. 

We want to choose $a$ and $u$ such that for any $s$ one would have 

\begin{equation}\label{beta}
1-\beta /3\leq |\Delta _1\cdots \Delta _s|\leq 1+\beta /3~~~
{\rm (recall~that~~}u=1+\beta {\rm )}~.
\end{equation}
As $\Delta _s=(\Delta _1\cdots \Delta _s)/(\Delta _1\cdots \Delta _{s-1})$, this  
implies  

\begin{equation}\label{beta/3}
1-\beta \leq \frac{1-\beta /3}{1+\beta /3}\leq 
|\Delta _s|\leq \frac{1+\beta /3}{1-\beta /3}\leq 1+\beta ~.
\end{equation}
The leftmost and rightmost inequalities are true for $\beta \in (0,1)$. 
Conditions (\ref{beta}) are fulfilled if 

\begin{equation}\label{ua}
ua/(1-a)+(ua)^2/(1-ua)\leq (u-1)/3
\end{equation}
which is true (together with conditions (\ref{twoconditions})) 
for $a=0.108$ and $u=1.7882$.

Now we can finish the proof of the convergence. 
Recall that the equations $(E_{j,i})$ 
were defined in the previous section. We represented the variables 
$\Delta _j$ as FPS in $q$ by iterating infinitely many times the 
following operation: a variable  
$\Delta _s$ is represented in the form 
$1-qT_{s,r}(q,\Delta ^r)$ 
using the equation $(E_{s,r})$ 
and then $1-qT_{s,r}(q,\Delta ^r)$ is substituted for 
$\Delta _s$ in all 
other equations of the infinite system. 
At any step we suppose that $|q|\leq a$ 
and $|\Delta _j|\in [1-\beta ,1+\beta ]$. The inequalities given 
in equation (\ref{beta/3}) 
imply that $|1-qT_{s,r}(q,\Delta ^r)|\in [1-\beta ,1+\beta ]$. 
We finally conclude that the series converge for $|q|\leq a$ and for all 
such $q$ one has $|\Delta _j|\in [1-\beta ,1+\beta ]$. 

For $a=0.108$, $\beta =0.7882$ one has 
$$|\Delta _{j+1}q^{j+1}|\leq (1+\beta )a|q^j|=(1.7882)(0.108)|q|^j
<(0.2118)|q|^j=(1-\beta )|q|^j\leq |\Delta _jq^j|$$ 
which implies that all zeros of $\theta (q,.)$ are distinct.~~~~~$\Box$

\section{Proof of Lemmas~\protect\ref{lemma1} and 
\protect\ref{lemma3}}

{\em Proof of Lemma~\ref{lemma1}:}\\ 

The second line of the formula is evident. Unless $\nu \leq \mu \leq s$ the 
claim of the lemma is trivial. 
It is also clear that 

$$(L_s^{-1})_{\mu ,\nu}=((-N_s)^{\mu -\nu})_{\mu ,\nu}=
(-1)^{\mu -\nu}
\Delta _s^{\mu -\nu}((N_s/\Delta _s)^{\mu -\nu})_{\mu ,\nu}~.$$ 
Set $P_s:=N_s/\Delta _s$. There remains to be proved that 

$$((P_s)^{\mu -\nu})_{\mu ,\nu}=
q^{(\mu -\nu)(s-\mu +1)+(\mu -\nu)(\mu -\nu -1)/2}~~
{\rm for}~~\nu \leq \mu \leq s~.$$ 
For $\mu -\nu =1$ this follows from the definition of $N_s$. Suppose that 
the above equation holds for $\mu -\nu \leq \kappa$. Note that 

$$((P_s)^{\mu -\nu})_{\mu ,\nu}=(P_s)_{\mu ,\mu -1}((P_s)^{\mu -\nu -1})_{\mu -1,\nu}~.$$
By induction the right-hand side equals

$$q^{s-\mu +1}q^{(\mu -\nu -1)(s-\mu +2)+(\mu -\nu -1)(\mu -\nu -2)/2}=
q^{(\mu -\nu)(s-\mu +1)+(\mu -\nu)(\mu -\nu -1)/2}$$
and the proof follows.~~~~$\Box$\\   

{\em Proof of Lemma~\ref{lemma3}:}\\ 

The coefficient $b_1$ equals $1+a+a^2+a^3+\cdots =1/(1-a)$ 
(independent of $s$). 
Suppose that the lemma is proved for $j\leq j_0<s-1$. Set 
$\Pi _1(a, M):=\prod _{k=s+1}^{\infty}(I+\sum _{j=1}^{s-1}(a^{k-s}M)^j)$ 
and present 
$\Pi _1$ in the form $I+c_1M+c_2M^2+\cdots +c_{s-1}M^{s-1}$. 
It is clear that $\Pi _1(a, M)=\Pi (a, aM)$. Therefore 
$c_j=b_ja^j$. In what follows we set $b_0=c_0=1$. Hence for $j\leq j_0$ one has  
$c_j\leq a^j/(1-a)(1-a^2)\cdots (1-a^j)$. Since  

$$\Pi (a,M)=(I+M+\cdots +M^{s-1})(I+c_1M+c_2M^2+\cdots +c_{s-1}M^{s-1})~,$$ 
one can obtain a term $M^{j_0+1}$ in one of the following ways:

1) one multiplies a term $M^{j_0+1}$ from one of the factors 
$I+\sum _{j=1}^{s-1}(a^{k-s}M)^j$ ($k\geq s$) by the terms $I$ of all the others; 
the sum of all these coefficients equals 
$\sum _{\nu =0}^{\infty}a^{\nu (j_0+1)}=1/(1-a^{j_0+1})$;

2) one multiplies $M^i$ from the factor $I+M+\cdots +M^{s-1}$ 
by $c_{j_0+1-i}M^{j_0+1-i}$ from 
the second factor for $i=1,\ldots ,j_0$. 
Thus 

$$\begin{array}{ccl}b_{j_0+1}&=&1/(1-a^{j_0+1})+c_1+\cdots +c_{j_0}\\ \\   
&\leq&
1/(1-a^{j_0+1})+\sum _{j=1}^{j_0}a^j/(1-a)(1-a^2)\cdots (1-a^j)\\ \\ 
&=&1/(1-a^{j_0+1})+\sum _{j=1}^{j_0}(1/(1-a)(1-a^2)\cdots (1-a^j)-
1/(1-a)(1-a^2)\cdots (1-a^{j-1}))\\ \\ 
&=&1/(1-a^{j_0+1})-1+1/B\\ \\
&=&a^{j_0+1}/(1-a^{j_0+1})+1/B~,\end{array}$$
where $B=(1-a)(1-a^2)\cdots (1-a^{j_0})$. One has 

$$a^{j_0+1}/(1-a^{j_0+1})+1/B=
(a^{j_0+1}(B-1)+1)/(1-a^{j_0+1})B<1/(1-a^{j_0+1})B$$
because $B\in (0,1)$. 
This proves the lemma by induction on $j$.~~~~~$\Box$\\

\end{document}